\documentclass[OJMO,manuscript]{cedram}

\usepackage[all]{xypic}

\usepackage{graphicx}%
\usepackage{multirow}%
\usepackage{amsmath,amssymb,amsfonts}%
\usepackage{amsthm}%
\usepackage{mathrsfs}%
\usepackage[title]{appendix}%
\usepackage{xcolor}%
\usepackage{textcomp}%
\usepackage{manyfoot}%
\usepackage{booktabs}%
\usepackage{algorithm}%
\usepackage{algorithmicx}%
\usepackage{algpseudocode}%
\usepackage{listings}%
\usepackage{caption}
\theoremstyle{thmstyleone}%
\theoremstyle{thmstyletwo}%

\theoremstyle{thmstylethree}%
\newcommand{\Null}{\mathrm{Null}}
\newcommand{\R}{\mathrm{R}}
\title[Article Title]{A New Active Set Scheme for Quadratic Programing}

\author{\firstname{Negin} \lastname{Bagherpour}}\address{Department of Engineering Sciences, University of Tehran, Tehran, 1417935840, Iran}\email{negin.bagherpour@ut.ac.ir}

\author{\firstname{Nima} \lastname{Minayi}}\address{Department of Engineering Sciences, University of Tehran, Tehran, 1417935840, Iran}\email{nima.minayi@ut.ac.ir}

\author{\firstname{AmirHossein} \lastname{Shanaghi}}\address{Department of Engineering Sciences, University of Tehran, Tehran, 1417935840, Iran}\email{shanaghi@ut.ac.ir}

  
\begin{abstract} 
We are faced with convex quadratic programing in many contexts related to control theory, economy and robotics. In this paper, we introduce a new active set algorithm for solving such problems and analyze its possible advantages. The novelty of the proposed scheme is in the way of solving the KKT system based on matrix properties. More precisely, we combine the two KKT equations to reduce the order and substitute it with a null space computation. The null space is in hand by use of the singular values decomposition. In problems with high number of independent constraints, we proposed another scheme. This also aims to solve the KKT system based on matrix properties. We implement both algorithms and test them over both randomly generated problems and standard problems mentioned in CUTEst. In general, more than 2000 tests with great variety are generated and computing times and accuracies are reported. The proposed schemes for solving convex quadratic problems are members of active set family. Because of using matrix properties it reduces computing time and as larger as the problem is, the improvement shows to be more remarkable. The first strategy performs the original acitive set when the number of constraints is low while the second outperforms the original algorithm when there exists a lot of independent constraints.
\end{abstract}

\begin{document}

\maketitle

\section{Introduction}\label{sec1}

Convex quadratic programming (CQP) is one of the most useful optimization problems. A lot of mathematical modelings such as mass inertia matrix, financial modeling and intelligent control are equivalent to a CQP. Moreover, sequential quadratic programing (SQP) is an iterative method for convex nonlinear programming which is based on quadratic approximation for the objective function and linear approximations for the constraints. So, in each iteration of SQP we in fact face with a CQP.

For a CQP with equality constraints, it is sufficient to solve the linear system resulting from the KKT conditions. However, for a general CQP active set algorithm shows to be effective \cite{R2}.  Different types of simplex are also advised for solving CQPs. Wolf's method \cite{R5}, Swarups simplex method \cite{R6} and Gupta and Sharma's method \cite{R7} are the most well known in this group.

The remainder of this paper is organized as follows: in Section \ref{sec2}, the original active set algorithm is reviewed. Then, in Section \ref{sec3} and Section \ref{sec4} our ideas for optimally solving the conditions of active set algorithm are stated. Finally, in Section \ref{sec5}, numerical results are reported to certify the efficiency of the proposed ideas for solving the CQPs faster than original active set algorithm. Our contributions are as follows: 1) We first review the original active set algorithm for solving CQPs. 2) We outline two ideas for solving the conditions resulting from KKT optimally. 3) Based on the reported numerical results both of the outlined ideas help to solve the KKT conditions faster. 

\section{Active Set Method}\label{sec2}
Lets consider the convex quadratic programing problem

\begin{eqnarray}\label{eqn1}
&\min& \frac{1}{2} X^T Q X+q^T X\nonumber\\
& s.t.& AX=b\nonumber\\
&& GX\leq h\nonumber
\end{eqnarray}

where $Q \in {\mathrm{R}}^{n\times n}$ is positive definite and $A \in {\mathrm{R}}^{m\times n}$, $G \in {\mathrm{R}}^{k\times n}$, $q \in {\mathrm{R}}^{n\times 1}$, $b \in {\mathrm{R}}^{m\times 1}$ and $h \in {\mathrm{R}}^{k\times 1}$ are arbitrary matrices and vectors. In the so-called active set algorithm, each iteration consists of the following main steps:\\

1) determining the current active set\\

2) computing the desired direction\\

3) computing the Lagrange multipliers and updating the active set\\

More precisely, with $X_0\in {\mathrm{R}}^n$ being the starting point, the current active set would be a matrix $A_0$ consisting of $A$ and $G(I)$ which refers to the rows $I$ of $G$ and satisfies $G(I)X_0=h(I)$. Hence, we have
\begin{eqnarray}\nonumber
A_0=\left[\begin{array}{c}
A\\
G(I)
\end{array}\right],\\
A_0X_0=\left[\begin{array}{c}
b\\
h(I)
\end{array}\right].
\end{eqnarray}
Now, a feasible point with enough reduction to $\Phi(X)=\frac{1}{2} X^T Q X+q^T X$ is needed. To compute such a direction, we solve the quadratic programing problem
\begin{eqnarray}\label{eqn2}
&\min& \frac{1}{2} P^TQP+X_0^TQP+q^TP.\\
&s.t.& A_0P=0\nonumber
\end{eqnarray}
The corresponding KKT system for (\ref{eqn2}) would be

\begin{eqnarray}
QP+A_0^T\lambda&=&-(QX_0+q)\label{eqn31}\\
A_0P&=&0\label{eqn32}
\end{eqnarray}

Before starting the next iteration, we need to update both the solution and the active set. The indices for which the Lagrange multipliers are positive will remain in active set while the constraint corresponding to the most negative Lagrange multiplier would be substitute by a new active constraint. In addition, the solution approximation would become $X=X_0+\alpha P$ where $\alpha$ is set to be the largest value for which $X$ is feasible. These steps are summarized in Algorithm \ref{alg1} \cite{R1}.

\bigskip

\begin{algorithm}
\caption{Active Set Algorithm for Convex Quadratic Programing}\label{alg1}
\begin{algorithmic}[1]
\State Compute the direction $P$ by solving (\ref{eqn2})
\If{$P=0$}
\State Compute the Lagrange multipliers $\lambda$ from (\ref{eqn31})
\If{$\lambda \geq 0$}
\State $X$ is the solution: \textbf{STOP}
\Else
\State remove the constraint corresponding to the smallest Lagrange multiplier
\EndIf
\Else
\State compute the maximum possible $\alpha$ value for which $X+\alpha P$ is a feasible solution
\State add the blocking constraint to the active set
\EndIf
\end{algorithmic}
\end{algorithm}

\section{New Scheme}\label{sec3}
The KKT conditions for solving convex quadratic problems might be established in the matrix form \cite{R2} to optimize the computations based on special matrix structures. In this paper, we develop two ideas to benefit from these special structures.We note that solving the QP subproblem (\ref{eqn2}) in Algorithm \ref{alg1} needs to be reasonably easier than the original QP problem. Here, a cost-effective algorithm is suggested to solve (\ref{eqn2}). Let rank of $A_0$ be equal to $k$ and its SVD decomposittion be $U\Sigma V_k^T$ with $U\in {\mathbb{R}}^{s\times s}$ and $V\in {\mathbb{R}}^{n\times n}$ being orthonormal and $\Sigma \in {\mathbb{R}}^{k\times k}$ being diagonal. We note that (\ref{eqn32}) gives $P=V_{n-k}Z\in \Null(A_0)$, while (\ref{eqn31}) results in 
\[P+r_0=V_kS\in \R(A_0^T),\]
where, $r_0=QX_0+q$. Hence, 
\[P=V_{n-k}Z=V_kS-r_0\]
and
\[
\left( V_k\hspace{0.3cm} V_{n-k}\right)
\left(\begin{array}{c}
S\\
-Z
\end{array}\right)
=r_0.\]
Since $V$ is orthonormal, we now get
\begin{equation}\label{eqn44}
P=-V_{n-k}V_{n-k}^Tr_0,
\end{equation}
which is the solution of (\ref{eqn2}). 

\section{Another Approach}\label{sec4}
We note that by use of a change of variable $\tilde{X}=Q^{\frac{1}{2}}X$, the coefficient matrix $Q$ can be eliminated in (\ref{eqn1}). So, without loss of generality, here we assume $Q=I$. The KKT condition for (\ref{eqn2}) will be 
\begin{eqnarray}
P+A_0^T\lambda&=&-(X_0+q)\label{eqn51}\\
A_0P&=&0\label{eqn52}
\end{eqnarray}
Substituting (\ref{eqn52}) in (\ref{eqn51}), we have $P^TP=-r_0^TP$ in which $r_0=X_0+q$. We then let $P=V_{n-k}Z$ and $\tilde{r_0}=V_{n-k}r_0$ to get
\begin{equation}\label{eqn6}
Z^TZ=-\tilde{r_0}^TZ.
\end{equation}
We note that the most important benefit of this strategy is reducing the dimension of the problem to $n-k$. So, this is of more interest in problems with greater number of independent constraints. The solution set of (\ref{eqn6}) clearly consists of the points on the sphere
\begin{equation}\label{eqn7}
\|Z+C\|=\|C\|
\end{equation}
where $C=\frac{\tilde{r_0}}{2}$. When $n-k$ is small enough, this approach is of significant interest. In case of $n-k=1$ and $n-k=2$ the results are given in Table \ref{tbl1}.
\begin{table}[ht]
\caption{Simple and common cases}\label{tbl1}
\begin{tabular}{c c c}\hline
 & $n-k=1$ & $n-k=2$\\\hline
 $Z$ & $Z=-2C$ & $Z=-C\pm \|C\|\left(\begin{array}{c}
\cos(\theta)\\
\sin(\theta)
\end{array}\right)$ \\
 $P$ & $P=V_{n-k}Z$ & $P=V_{n-k}Z$\\
\hline\end{tabular}
\end{table}
Let describe the process in details. First, the argument $\theta$ needs to be computed so that
\[V_{n-k}^Tr_0=-Z=C\mp \|C\|\left(\begin{array}{c}
\cos(\theta)\\
\sin(\theta)
\end{array}\right).\]
Then, $P$ is calculated from $Z$. So far, we presented two schemes for solving the subproblem presented in Algorithm \ref{alg1} for computing the direction $P$. Theses are summarized in Table \ref{tbl2} below.
\begin{table}[ht]
\caption{Strategies for computing $P$}\label{tbl2}
\begin{tabular}{c c c }\hline
& Schme 1 (Section \ref{sec3})& Scheme 2 (Section \ref{sec4})\\\hline
usage & low $n-k$ & otherwise\\
$P$ & $P=V_{n-k}Z$ with $Z$ on a sphere & $P=-V_{n-k}V_{n-k}^Tr_0$\\
\hline
\end{tabular}
\end{table}

Next, in Section \ref{sec5} numerical test results are reported to confirm the efficiency of Algorithm \ref{alg1} by both schemes for computing $P$.

\section{Numerical Results}\label{sec5}
\begin{figure}[h!]
\centering
\begin{subfigure}{0.43\textwidth}
    \includegraphics[width=\textwidth]{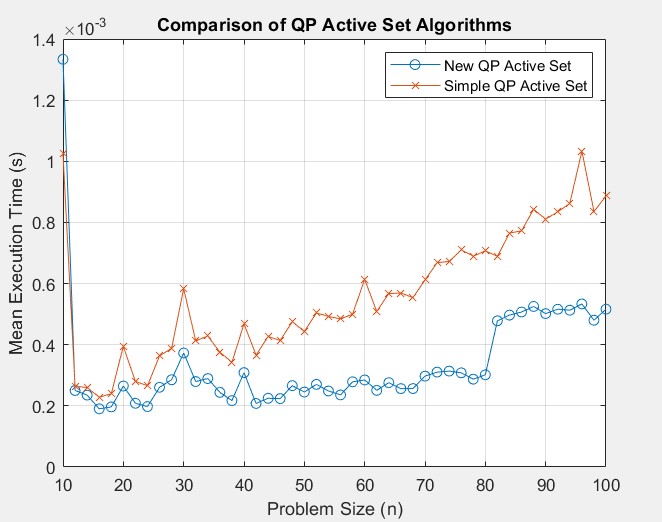}
    \caption{$10\leq n \leq 100$, $n_e=1$ and $n_i=10$}
    \label{fig:1}
\end{subfigure}
\hfill
\begin{subfigure}{0.43\textwidth}
    \includegraphics[width=\textwidth]{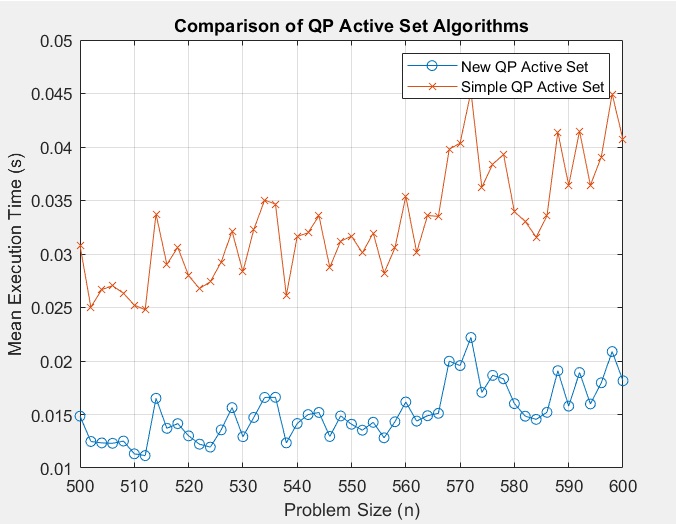}
    \caption{$500\leq n \leq 600$, $n_e=1$ and $n_i=10$}
    \label{fig:2}
\end{subfigure}
\begin{subfigure}{0.43\textwidth}
    \includegraphics[width=\textwidth]{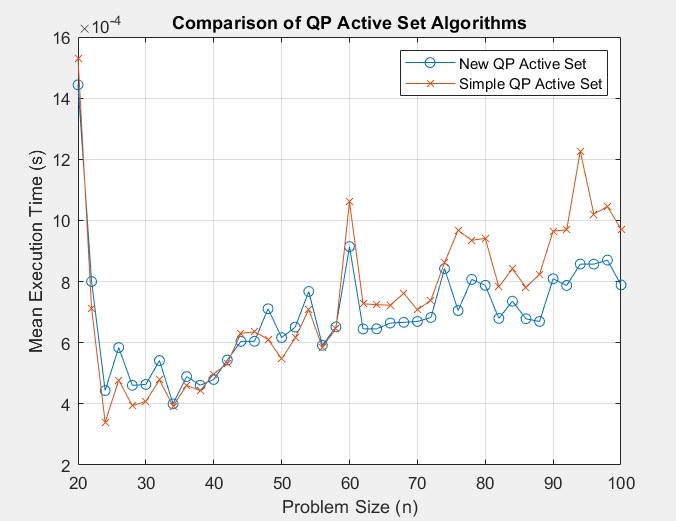}
    \caption{$20\leq n \leq 100$, $n_e=10$ and $n_i=10$}
    \label{fig:3}
\end{subfigure}
\hfill
\begin{subfigure}{0.43\textwidth}
    \includegraphics[width=\textwidth]{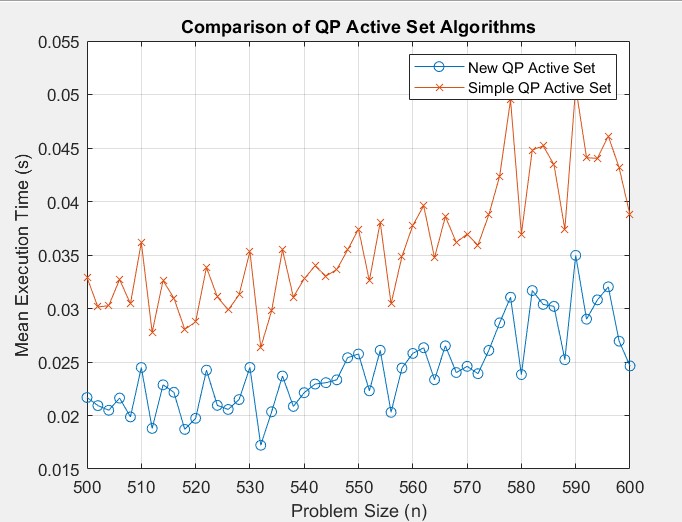}
    \caption{$500\leq n \leq 600$, $n_e=10$ and $n_i=10$}
    \label{fig:4}
\end{subfigure}        
\begin{subfigure}{0.43\textwidth}
    \includegraphics[width=\textwidth]{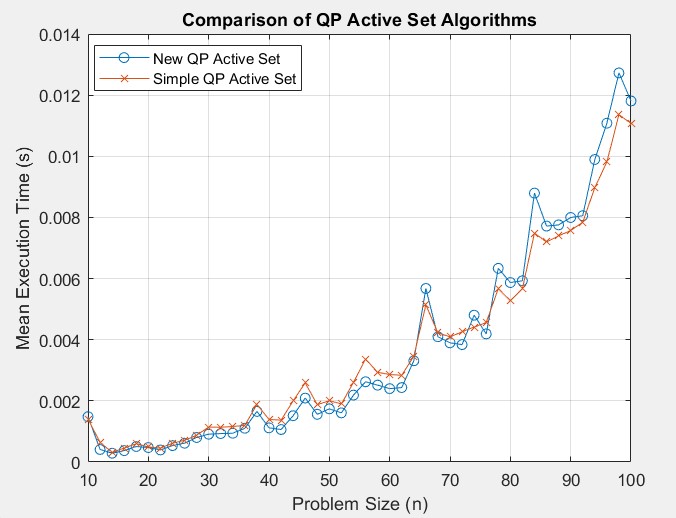}
    \small\caption{$10\leq n \leq 100$, $n_e=n-1$ and $n_i=1$}
    \label{fig:5}
\end{subfigure}
\hfill
\begin{subfigure}{0.43\textwidth}
    \includegraphics[width=\textwidth]{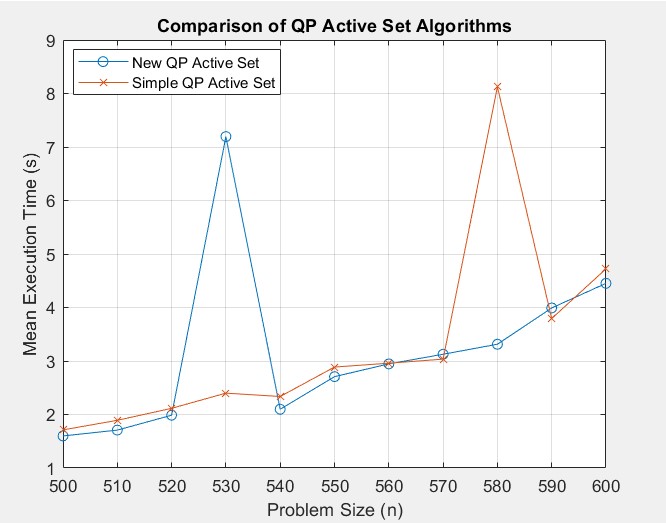}
    \caption{$100\leq n \leq 600$, $n_e=n-1$ and $n_i=1$}
    \label{fig:6}
\end{subfigure}
\begin{subfigure}{0.43\textwidth}
    \includegraphics[width=\textwidth]{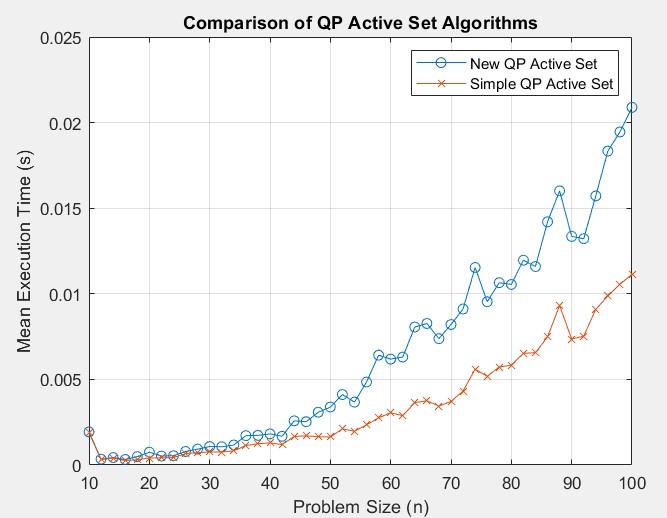}
    \caption{$10\leq n \leq 100$, $n_e=\frac{n}{2}$ and $n_i=\frac{n}{2}$}
    \label{fig:7}
\end{subfigure}
\hfill
\begin{subfigure}{0.43\textwidth}
    \includegraphics[width=\textwidth]{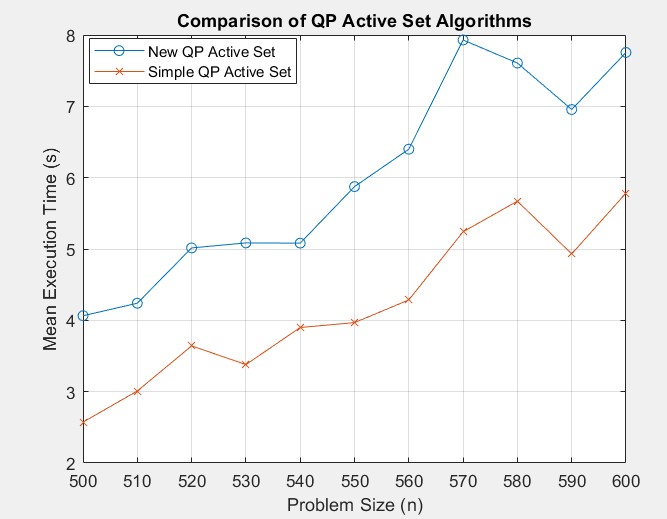}
    \caption{$500\leq n \leq 600$, $n_e=\frac{n}{2}$ and $n_i=\frac{n}{2}$}
    \label{fig:8}
\end{subfigure}        
\caption{Time Comparison of the Original Active Set and Scheme 1}
\label{fig:figures1}
\end{figure}

We first used Scheme 1 to solve 1000 randomly generated test problems with up to 600 unknown variables. The random tests are generated with a good variety with different number of equality and inequality constraints, $n_e$ and $n_i$. The computing times for the original active set algorithm is compared with Scheme 1 in Figure \ref{fig:figures1}.

\begin{figure}[ht!]
\centering
\begin{subfigure}{0.43\textwidth}
    \includegraphics[width=0.88\textwidth]{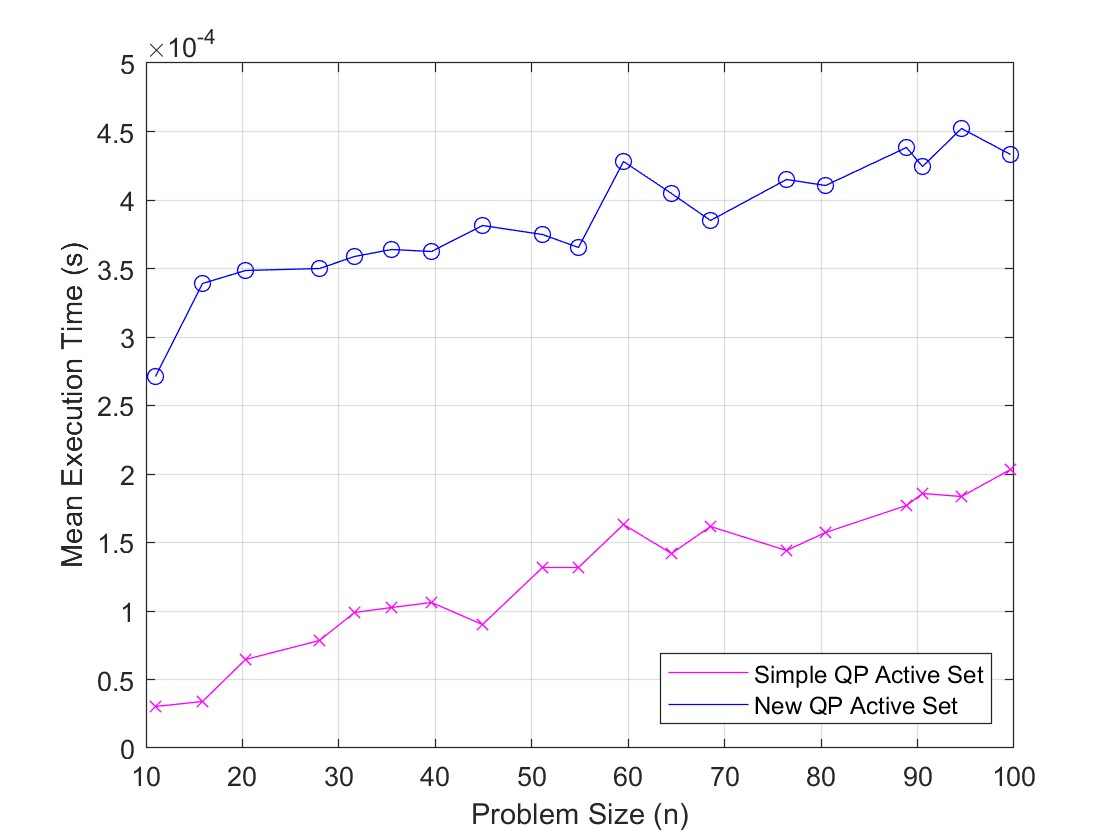}
    \caption{$10\leq n \leq 100$, $n_e=1$ and $n_i=10$}
    \label{fig:21}
\end{subfigure}
\hfill
\begin{subfigure}{0.43\textwidth}
    \includegraphics[width=0.88\textwidth]{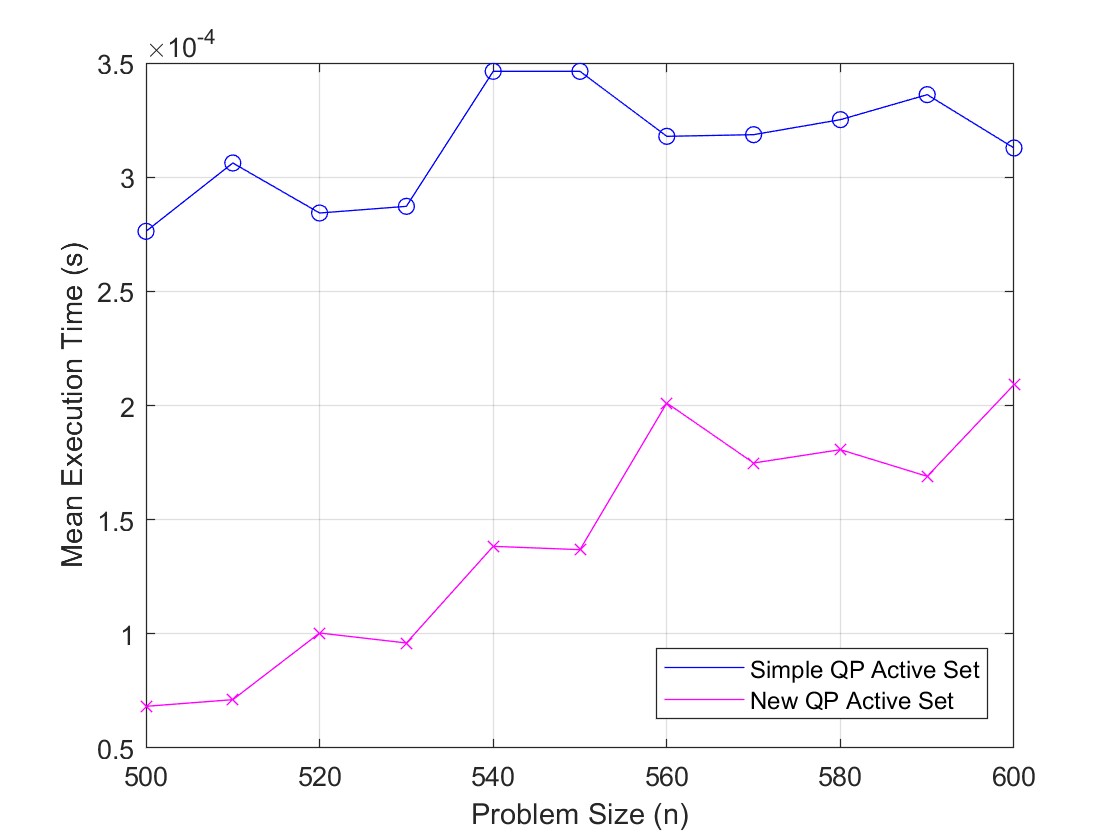}
    \caption{$500\leq n \leq 600$, $n_e=1$ and $n_i=10$}
    \label{fig:22}
\end{subfigure}
\begin{subfigure}{0.43\textwidth}
    \includegraphics[width=0.88\textwidth]{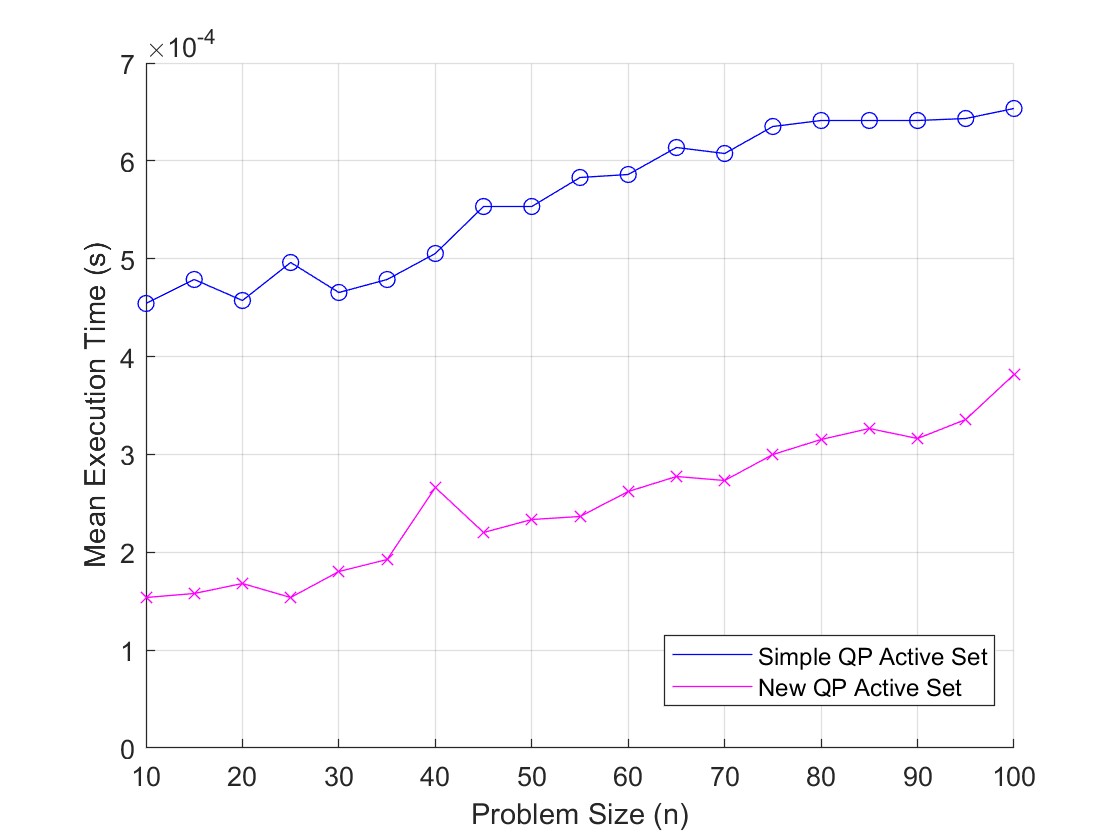}
    \caption{$20\leq n \leq 100$, $n_e=10$ and $n_i=10$}
    \label{fig:23}
\end{subfigure}
\hfill
\begin{subfigure}{0.43\textwidth}
    \includegraphics[width=0.88\textwidth]{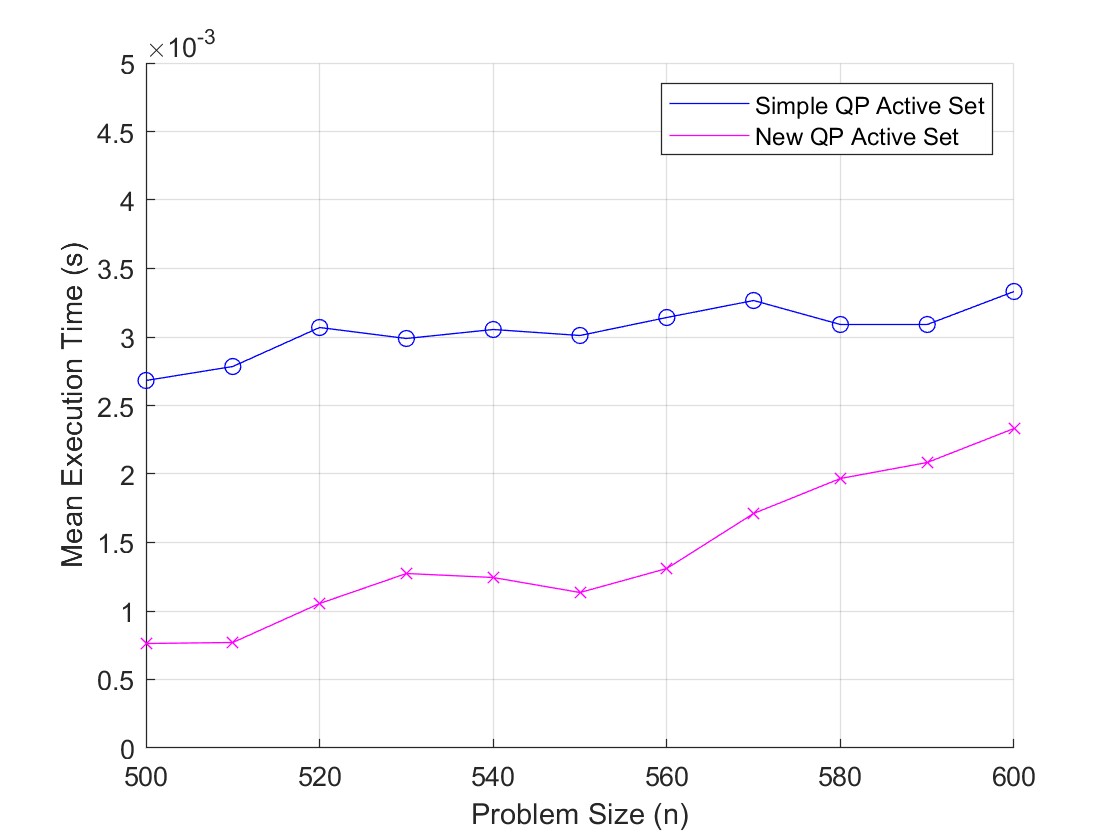}
    \caption{$500\leq n \leq 600$, $n_e=10$ and $n_i=10$}
    \label{fig:24}
\end{subfigure}        
\begin{subfigure}{0.43\textwidth}
    \includegraphics[width=0.88\textwidth]{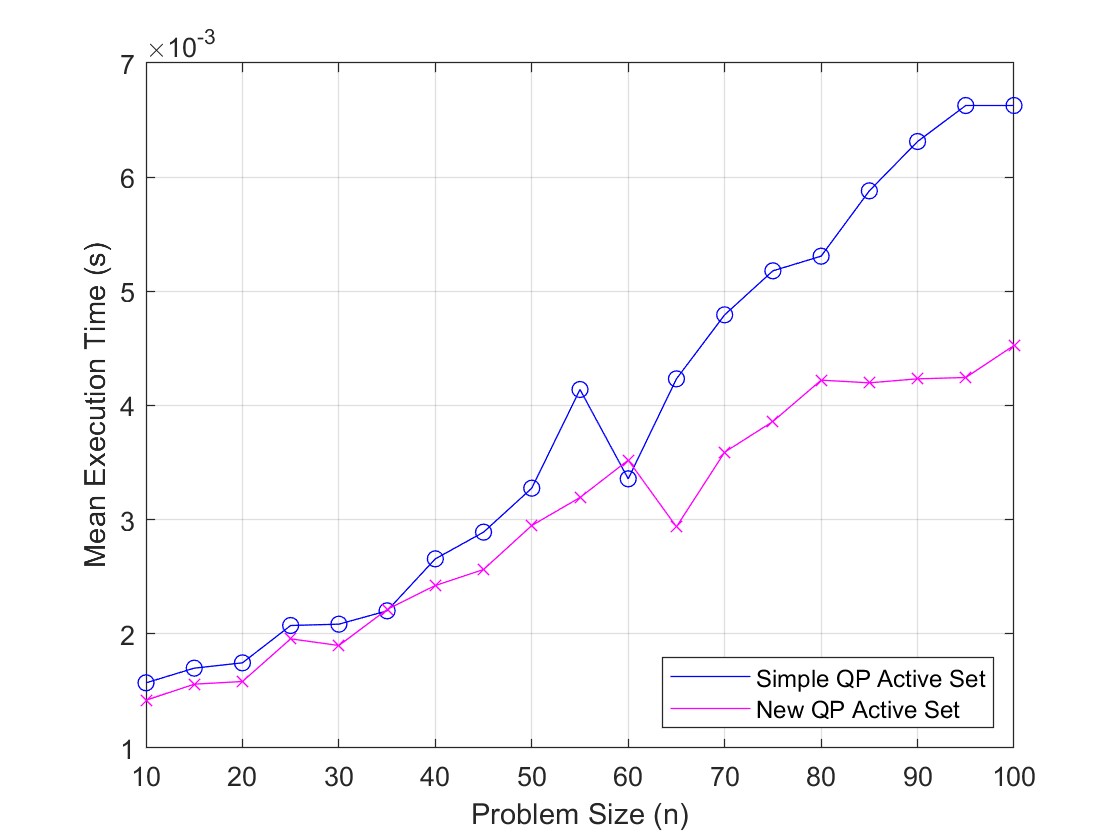}
    \caption{$10\leq n \leq 100$, $n_e=n-1$ and $n_i=1$}
    \label{fig:25}
\end{subfigure}
\hfill
\begin{subfigure}{0.43\textwidth}
    \includegraphics[width=0.88\textwidth]{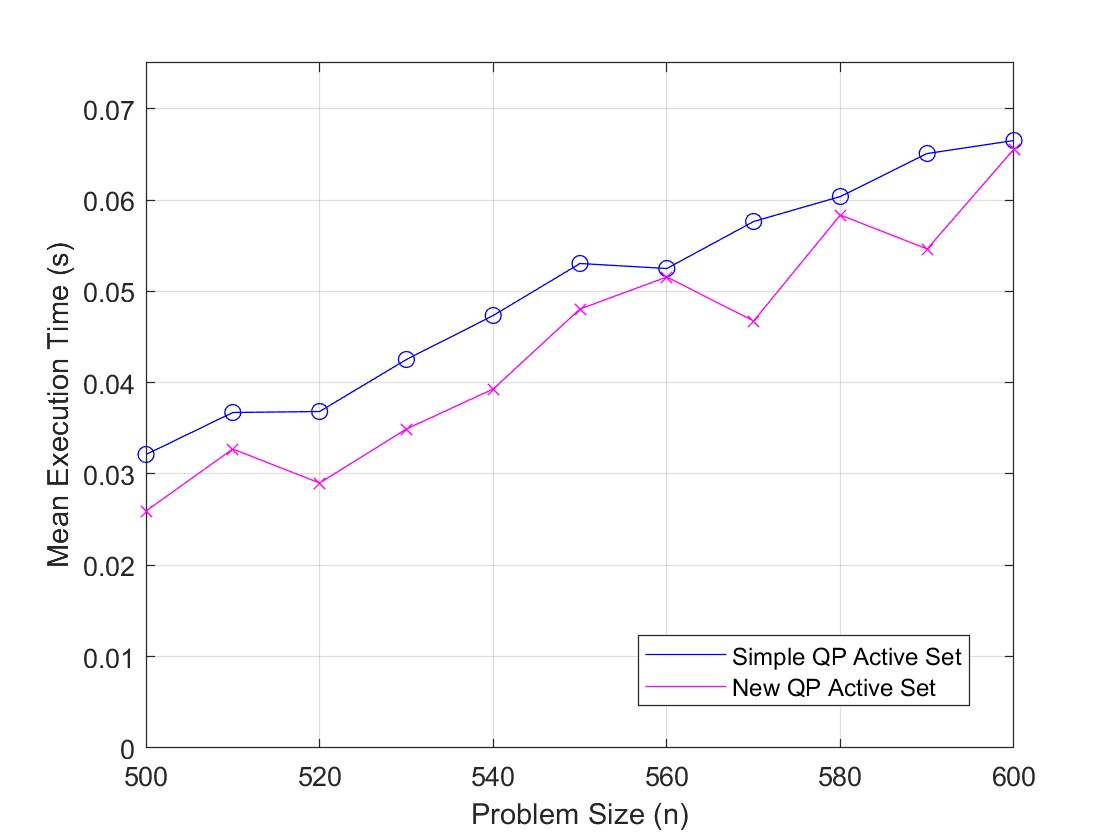}
    \caption{$100\leq n \leq 600$, $n_e=n-1$ and $n_i=1$}
    \label{fig:26}
\end{subfigure}
\begin{subfigure}{0.43\textwidth}
    \includegraphics[width=0.88\textwidth]{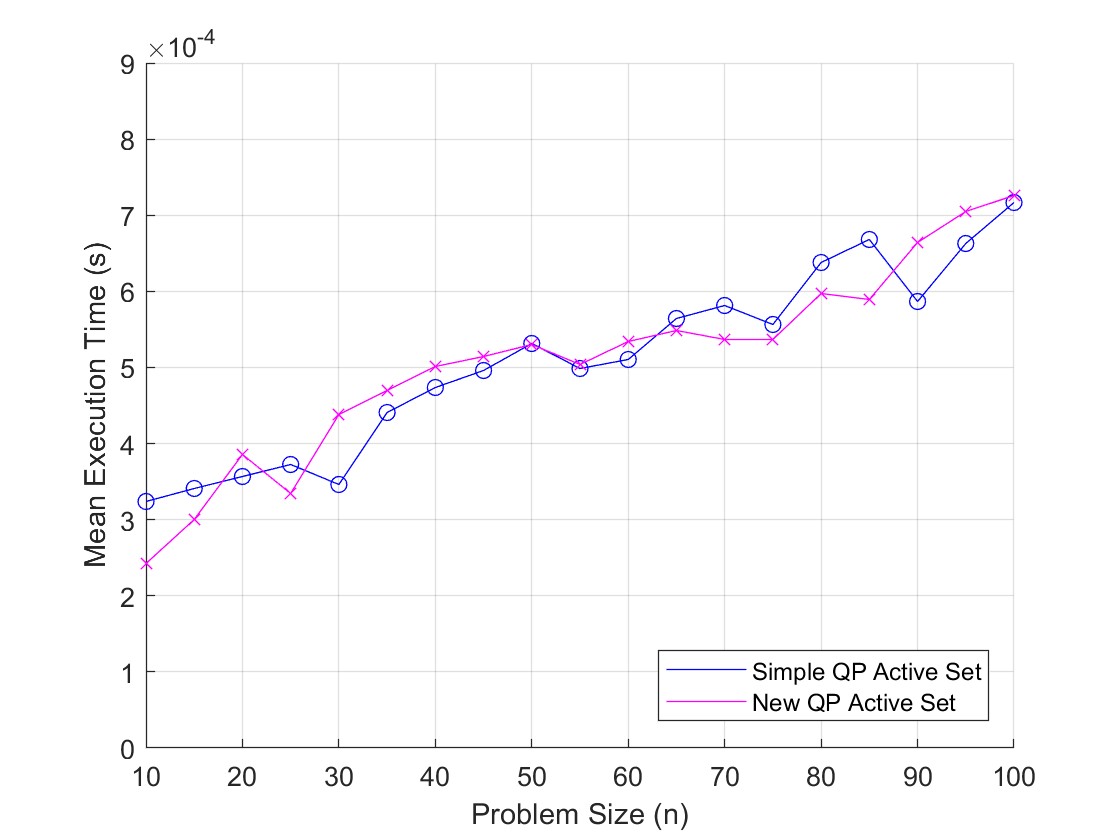}
    \caption{$10\leq n \leq 100$, $n_e=\frac{n}{2}$ and $n_i=\frac{n}{2}$}
    \label{fig:27}
\end{subfigure}
\hfill
\begin{subfigure}{0.43\textwidth}
    \includegraphics[width=0.88\textwidth]{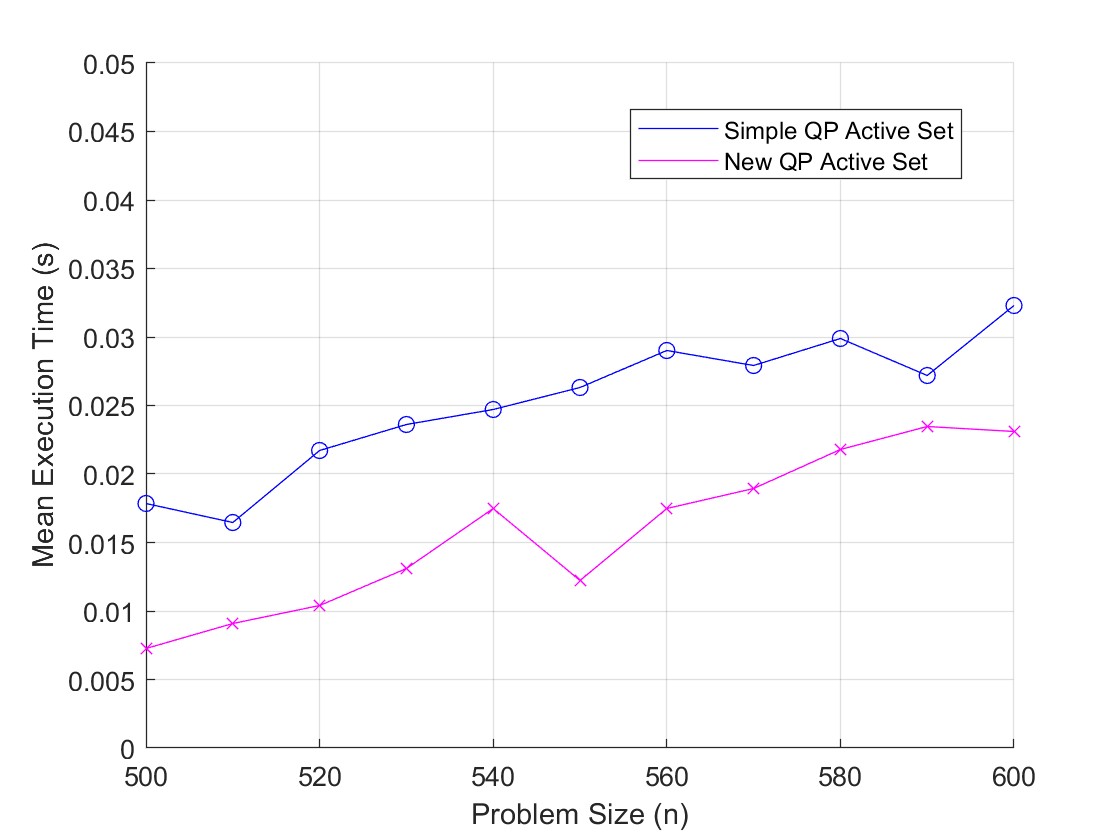}
    \caption{$500\leq n \leq 600$, $n_e=\frac{n}{2}$ and $n_i=\frac{n}{2}$}
    \label{fig:28}
\end{subfigure}        
\caption{Time Comparison of the Original Active Set and Scheme 2}
\label{fig:figures2}
\end{figure}

For larger values of $n$ and small number of constraints, Scheme 1 compute the solution faster. According to figures \ref{fig:7} and \ref{fig:8}, Scheme 1 fails to compute the solution in a competetive time. Now, we compare the computing times of the original active set and Scheme 2 in Figure \ref{fig:figures2}. In this experiment, tests are conducted with high number of independent constraints.

\vfill
\newpage
\clearpage

\begin{figure}[ht!]
\centering
\includegraphics[height=4.5cm]{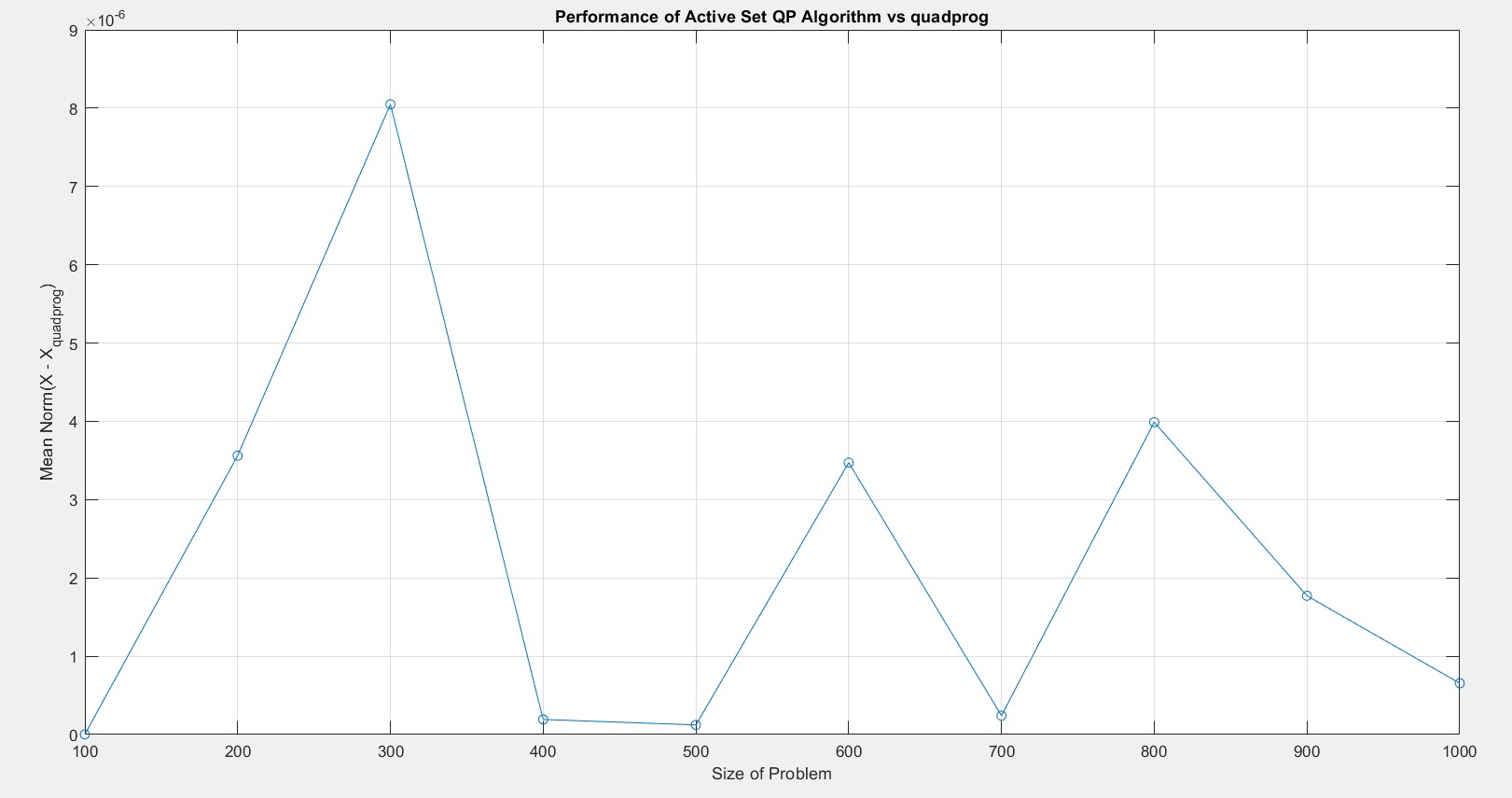}
\caption{Error Norm for Scheme 1}\label{fig39}
\end{figure}

\newpage
\clearpage
Since, Scheme 1 failed to solve some of the test problems in the expected time, we need to confirm its convergence regardless of computing times. For the proposed scheme, we present the average error norm in Figure \ref{fig39} to confirm that although the computing time is not competetive, the convergence has been achieved. The error norm here is defined to be the norm of difference between the computed solution and the exact solution computed by quadprog command of MATLAB.

\begin{figure}[ht!]
\centering
\includegraphics[height=4.5cm]{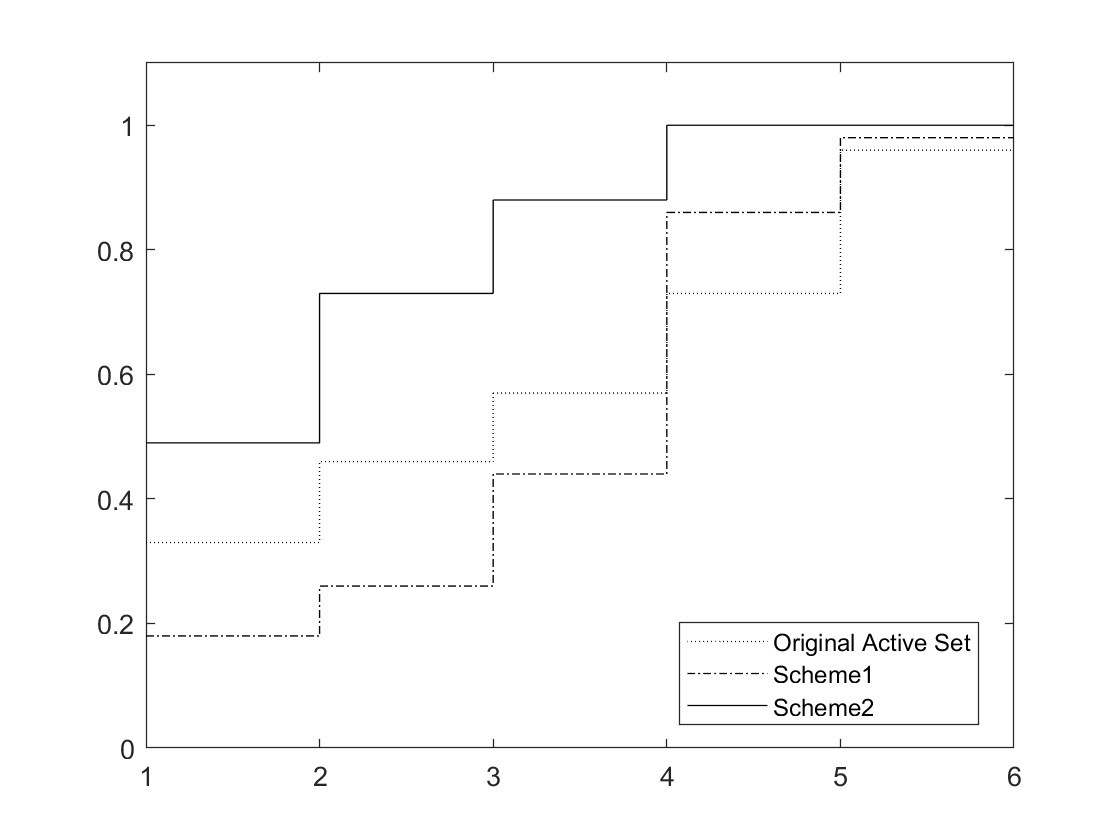}
\caption{Dolan More Time Comparison}\label{fig44}
\end{figure}

Finally, in Figure \ref{fig44}, the Dolan More \cite{R3} time profiles are presented to compare Algorithm \ref{alg1} with scheme 1 and 2 in solving the QP standard tests introduced in CUTEst \cite{R4}. Considering the rpesented profile, it can be seen that Sceme 2 outperforms other in converging to the solution faster. Moreover, it is the only method (between the three) which is able to solve all the CUTEst test problems in at most four time the best.

\bibliographystyle{plain+eid}

\nocite{*}

\bibliography{samplebib}

\begin{thebibliography}{1}

\bibitem{R3}
E.~Dolan and J.~B. More.
\newblock Benchmarking optimization software with performance profiles.
\newblock {\em Math. Prog.}, 91:201--213, 2002.

\bibitem{R4}
N.~I.~M. Gould, D.~Orban, and Ph. Toint.
\newblock Cutest: A constrained and unconstrained testing environment, latest,
  Jun 2022.

\bibitem{R7}
A.~Gupta and J.~Sharma.
\newblock A generalized simplex technique for solving quadratic-programming
  problem.
\newblock {\em Indian J. Technol.}, 21:198--201, 1983.

\bibitem{R2}
Z.~Iqbal and S.~Nooshabadi.
\newblock Exploiting block structures of kkt matrices for efficient solution of
  convex optimization problems.
\newblock {\em IEEE Access}, 9:116604--116611, 2021.

\bibitem{R1}
J.~Nocedal and S.~Wright.
\newblock {\em Numerical Optimization}.
\newblock Springer, New York, 2006.

\bibitem{R6}
N.~Suleiman and M.~Nawkhass.
\newblock A new modified simplex method to solve quadratic fractional
  programming problem and compared it to a traditional simplex method by using
  pseudoaffinity of quadratic fractional functions.
\newblock {\em Appl. Math. Sci.}, 7:3749--3764, 2013.

\bibitem{R5}
Ph. Wolfe.
\newblock The simplex method for quadratic programming.
\newblock {\em Econometrica}, 27:382--398, 1959.

\end{thebibliography}

\end{document}